\documentclass[12pt,a4paper,reqno]{amsart}

\usepackage[utf8]{inputenc}
\usepackage{color}
\usepackage{enumerate}
\usepackage[margin=2.5cm]{geometry}
\usepackage{graphicx} 
\usepackage{overpic}
\usepackage{amsaddr}
\usepackage{cite}

\usepackage[active]{srcltx}
\usepackage[colorlinks=true,breaklinks=true,linkcolor=blue,urlcolor=green,citecolor=blue]{hyperref}
\usepackage[utf8]{inputenc}
\usepackage[USenglish]{babel}
\usepackage{algorithm}
\usepackage{amsfonts}
\usepackage{amsmath,amssymb,amsthm,mathabx}
\usepackage{latexsym}
\usepackage{lineno}
\usepackage{lmodern}
\usepackage{mathrsfs}
\usepackage{multirow}
\usepackage{mwe}
\usepackage{pgf,tikz}
\usepackage{psfrag}
\usepackage{soul}
\usepackage{subcaption}
\usepackage{yfonts}
\usetikzlibrary{arrows}

\RequirePackage{lineno}

\input{apalike}

\usepackage[colorinlistoftodos]{todonotes}
\usepackage[normalem]{ulem}

\definecolor{ColorErnest}{rgb}{1.0,0,0.}
\definecolor{ColorToni}{rgb}{0.1,0.5,0.2}
\definecolor{ColorMT}{rgb}{0.9,0.2,0.1}
\definecolor{Bluish}{rgb}{0.,0.,0.5}
\definecolor{ChaosBlue}{HTML}{51abe3}
\definecolor{grayish}{rgb}{0.5,0.5,0.5}
\definecolor{Reddish}{rgb}{0.5,0.,0.}
\definecolor{whitish}{rgb}{0.75,0.75,0.75}
\definecolor{ColorST}{rgb}{0.8,0.8,0.8}

\newcommand\tgcom[2][noinline]{\todo[#1, color=ColorToni!20!white]{\small \texttt{TG}: #2}}

\newcommand{\efcom}[2][noinline]{\todo[#1, color=ColorErnest!15!white]{\small \texttt{EF}: #2}}

\newcommand{\stkout}[1]{\ifmmode\text{\sout{\ensuremath{#1}}}\else\sout{#1}\fi}

\newtheorem{theorem}{Theorem}
\newtheorem{remark}{Remark}

\newtheorem{definition}{Definition}

\newcommand{\Id}{\mathrm{Id}}

\begin{document}

\title[Functional shifts in cooperative systems]{Functional shift-induced degenerate transcritical Neimark-Sacker bifurcation in a discrete hypercycle}

\author[E. Fontich, A. Guillamon, J. Perona, J. Sardanyés]
{{Ernest Fontich$^{1,2,**}$,  Antoni Guillamon$^{3,2}$, Júlia Perona$^{3}$, \MakeLowercase{and} Josep Sardanyés$^{2}$}}
\address{$^{1}$Departament de Matem\`atiques i Inform\`atica, Universitat de Barcelona (UB), Gran Via de les Corts Catalanes 585, 08007 Barcelona, Spain}
\thanks{$^{**}$Corresponding author}
\email{fontich@ub.edu}
\address{$^{2}$Centre de Recerca Matem\`atica (CRM), 
08193 Bellaterra, Barcelona, Catalonia, Spain }
\address{$^{3}$Departament de Matem\`atiques and IMTech, Universitat Polit\`ecnica de Catalunya, EPSEB, Av. Dr. Mara\~n\'on 44-50, 08028 Barcelona, Catalonia, Spain}


\maketitle

\section*{Abstract}
 In this article we investigate the impact of functional shifts in a time-discrete cross-catalytic system. We use the hypercycle model considering that one of the species shifts from a cooperator to a degradader. At the bifurcation caused by this functional shift, an invariant curve collapses to a point $P$ while, simultaneously, two fixed points collide with $P$ in a transcritical manner. All points of a line containing $P$ become fixed points at the bifurcation and only at the bifurcation. Hofbauer and Iooss~\cite{HofbauerIooss1984} presented and proved a result that provides sufficient conditions for a Neimark-Sacker bifurcation (the authors called it ``Hopf") to occur in a special degenerate situation. They use it to prove the existence of an invariant curve for the model when a parameter related to the time discreteness of the system goes to infinity becoming a continuous-time system. Here we study the bifurcation that governs the functional shift and demonstrate the existence of an invariant curve when the cooperation parameter approaches zero and thus approaches the switch to degrading species. This invariant curve lives in a different domain and exists for a different set of values of the parameters described by these authors. In order to apply the mentioned result we uncouple the Neimark-Sacker and the transcritical bifurcations. This is accomplished by a preliminary singular change of coordinates that puts the involved fixed points at a fixed position, so that they stay at a fixed distance among them. Finally, going back to the original variables, we can describe mathematically the details of this bifurcation.


\section*{Introduction}
Hypercycles are catalytic sets of macromolecules, where each replicator catalyzes the replication of the next species of the set. This concept was first introduced by Manfred Eigen and Peter Schuster in 1977~\cite{EigenManfred1977} and has played a pivotal role in the study of prebiotic evolution and the overcoming of the so-called information crisis~\cite{Eigen1971,Eigen1982,Smith1995}. Research in hypercycles primarily investigates cooperative interactions among replicators~\cite{Eigen1979}. Hypercycle theory has been also applied to investigate dynamics of ecological systems~\cite{Smith1995,Sardanyes2009} and has aided in modeling experimental systems undergoing cooperation in engineered bacteria~\cite{Amor2017}. Cooperation has been previously described in different experimental systems with coiled-coil peptides~\cite{Lee1997}, yeast cell populations~\cite{Shou2007}, and self-replicating ribozymes~\cite{Vaidya2002}. Despite cooperation is a major driver in both molecular and ecological interactions, additional interactions are likely to emerge. For instance, molecular catalytic replicators experience mutational processes that could change functional properties, while ecological species may undergo behavioral shifts in response to environmental and ecological changes. These shifts may entail transitions from cooperative to antagonistic interactions, which we refer to as \emph{functional shifts}. The dynamics arising from these functional shifts in catalytic cycles has been scarcely investigated. Few works have addressed this subject for continuous-~\cite{Boerlijst1995,Bassols2021} and discrete-time~\cite{Perona2020} systems. 

Examples of functional shifts are widespread in Ecology. In the marine environment, large pelagic predators like tuna, sharks, and dolphins are known to collaborate in locating and handling small pelagic shoals~\cite{Scott2012}. However, instances of predation among these species have also been documented~\cite{Maldini2003,Melillo2014}. Similarly, seabirds form inter-specific flocks to cooperate in finding food at sea, but they may also engage in predation during breeding in colonies~\cite{Votier2004,Almaraz2011}. Waterbirds, on the other hand, form mixed colonies to protect against predators, but certain species may opportunistically prey on other species, especially under adverse environmental conditions~\cite{Hatch1970,Andersson1976}. These examples highlight the complexity of ecological interactions and how cooperation and antagonism can coexist within diverse ecosystems. Despite the existence of such observations, the specific effects of functional shifts on the dynamics and stability of cooperating systems have not been thoroughly explored. 

Dynamics in hypercycles have been extensively studied by means of continuous-time mathematical models, both in the limit of infinite diffusion~\cite{Eigen1979,Campos2000,Sardanyes2006,Silvestre2008,Puig2018} and in spatially extended systems~\cite{Boerlijst1995,Cronhjort1995}. The use of discrete models has been less explored although these systems are of great importance for species with non-overlapping generations such as insects in boreal climates.  The first author to investigate a discrete model was Hofbauer~\cite{Hofbauer1984}. In his model, the interactions between species are denoted by $k_i$, being $x_i$ the concentration of each species. Since the work by Hofbauer, very few research has been performed in discrete-time hypercycles~\cite{Hofbauer1984,HofbauerIooss1984,Perona2020}. Specifically, in Ref.~\cite{Perona2020} we went in depth into Hofbauer's discrete model. We proved both the cooperation in a $n$-dimensional system (i.e., no species goes extinct) and the existence of a fixed point in the interior of the domain. This fixed point is {globally asymptotically stable} in the three-dimensional system and unstable for $n\geq 4$. We also carried out a numerical study to find invariant curves for the four-dimensional system. To compute the invariant curve we built a pseudo Poincar\'e map inspired by~\cite{Gelfreich2019} taking advantage that the discrete system is close to a continuous one. Our main contribution was the investigation of dynamics and the bifurcations when a species of the hypercycle shifts from cooperation, $k_i > 0$, to degradation $k_i < 0$, in small hypercycles i.e., $n = 2,3,4$. We assumed that the first species was the one suffering the shift, i.e., $k_1=0$. We analytically proved that the trajectories tend asymptotically to a fixed point in a corner of the domain when the first species shifts to degradation, for any number of species in the hypercycle. We analytically obtained the rates of convergence to the fixed points of the studied hypercycles. We also numerically obtained that the discrete-time hypercycle is governed by an invariant attracting curve that shrinks to a corner of the domain and disappears throughout a Neimark-Sacker bifurcation when $k_1= 0$. Moreover, the interior fixed point collides with the corner point in a transcritical bifurcation and a line of fixed points appears, thus making the bifurcation more degenerate.

The main goal of this paper is to analytically prove the existence of the aforementioned invariant curves for {the} discrete-time hypercycle introduced by Hofbauer \cite{Hofbauer1984} by inspecting the cooperation parameter that drives the functional shift. We first introduce the {model and, for} the sake of completeness, we compute the fixed points and their stability as a function of the parameters. Next, we recall the Neimark-Sacker bifurcation for discrete-time dynamical systems and a version of it due to Hofbauer and Iooss~\cite{HofbauerIooss1984} that proves the existence of a family of attracting invariant curves in a family of maps that can be expressed as a step of the Euler's integration method of a differential equation; the corresponding vector field has a fixed point with a pair of purely imaginary eigenvalues while the other eigenvalues have negative real part and, moreover, the real part of the coefficient of the resonant term of lowest degree is negative. Finally, we apply the theorem by Hofbauer and Iooss to our discrete-time hypercycle with four species, $n=4$, for $k_1\to 0^{+}$. To do so, we make a singular change of coordinates to make our system less degenerate. We also carry out a translation to place the fixed point at the origin and we rewrite the system in the form stated in the hypothesis of the theorem. Then, we prove both hypotheses of the theorem to conclude that the four-species system has an attracting invariant curve that appears when $k_1=0$ through a degenerate Neimark-Sacker bifurcation.

\section{Hofbauer's discrete-time hypercycle model}\label{sect:Hofbauer_discrete_model}

In this section we present the discrete-time model for the hypercycle, introduced by Hofbauer in \cite{Hofbauer1984}, and we relate it to a continuous-time model. We also compute the basic elements of the dynamics such as the fixed points of the system and their stability. This dynamical system consists of a set of $n$ species $s_i$, $1\le i\le n$, such that the species $s_{i-1}$ catalyzes {only} the next one $s_i$, in a cyclic manner, with a strength  $k_i$. Let
$x_i$ be the concentration of the $i$-th species.
For convenience of notation, we write $x_0:=x_n$ and $x_{n+1}:=x_1$ and similarly for $k_0,\, k_{n+1}$.
The model assumes that if the total population is normalized to 1, it remains constant. This is accomplished by the introduction of a flux
$\phi(x)$, which also introduces competition between all the hypercycle species. This fact implies that the system will be defined on the $n$-simplex
\begin{equation}\label{eq62}
	S_n=\{x=(x_1,\dots,x_n)\in\mathbb{R}^n | \; \sum_{i=1}^nx_i=1, \, x_i\geq 0, \, 1\le i\le n\}.
\end{equation}
We introduce the hyperplane
\begin{equation*}
	\Delta_n=\{x\in\mathbb{R}^n  | \; \sum_{i=1}^n x_i=1 \},
\end{equation*}
and the set $
	\widetilde{\Delta}_n=\{x\in\Delta_n  | \; x_i\neq 0, \,  1\le i\le n\}.
$	

The system is determined by the map $F=(F_1, \dots ,F_n): S_n\to S_n$, where
\begin{equation}\label{eq63}
	F_i(x)=\frac{C+k_ix_{i-1}}{C+\phi(x)}x_i, \qquad 1\le i\le n,
\end{equation}
$C>0$ is a constant of proportionality and
\begin{equation}\label{eq61p5}
	\phi(x)=\sum_{i=1}^nk_ix_ix_{i-1}.
\end{equation}
In \cite{Hofbauer1984} the map \eqref{eq63} is related to the corresponding 
continuous-time system
\begin{equation}\label{eq61}
	\dot{x_i}=x_i(k_ix_{i-1}-\phi(x)), \qquad 1\le i\le n,
\end{equation}
which also satisfies that if the initial total population is 1, then it remains constant.
In particular we can write
\begin{equation}\label{eq63p5}
	\dfrac{F_i(x)-x_i}{C^{-1}}=C\,\left(\frac{C+k_ix_{i-1}}{C+\phi(x)}x_i-x_i\right) =x_i\left(k_ix_{i-1}-\phi(x)\right)\frac{C}{C+\phi(x)}.
\end{equation}
This expression allows us to compare the map $F$ with the continuous time model, since $C^{-1}$ can be interpreted as the time interval between two generations and $F(x)$ can be seen as the Euler step of length $C^{-1}$ of the continuous-time system (\ref{eq61}) since
\begin{equation*}
	\lim_{C^{-1} \to 0}\frac{x_i(t+C^{-1})-x_i(t)}{C^{-1}}=x_i{(t)}(k_i\,x_{i-1}{(t)}-\phi(x(t))),
\end{equation*}
identifying $F_i(x)(t)$ with $x_i(t+C^{-1})$ provided that $x_i=x_i(t)$. Then, for large values of $C$, the discrete-time model (\ref{eq63}) approximates the differential equation (\ref{eq61}) and
therefore we expect that, in such case, both models have similar properties.

If we keep the constants $k_i$ positive and bounded away from zero, and we let one of them, say $k_\ell$, go to zero and then become negative, the model can be interpreted biologically as there has been a \emph{functional shift} meaning that the role of the species $s_{\ell -1}$ changes from cooperation  ($k_\ell >0$) to degradation ($k_\ell <0$). Note that, due to the cyclic character of our model, we can assume, without loss of generality, that the parameter that tends to zero is $k_1$ while all $k_i$ with $i>1$ are bounded away from zero.

In Ref.~\cite{Perona2020} this bifurcation was studied  and it was found numerically that for four species there is an attracting invariant curve that tends to the point $Q=(0,{0} ,0,1) $ when $k_1 \to 0^{+}$. Also, when $k_1>0$  there is a unique fixed point $P$ in  the interior of the simplex $S_n\cap \widetilde \Delta_n$, already described in \cite{Hofbauer1984}. The point $P$ collides with $Q$ when $k_1=0$, and goes out of $S_n$ in a transcritical-like bifurcation.
Moreover, in \cite{Perona2020} it is shown that for any number of species, when $k_1\le 0$, the basin of attraction of $Q$ contains $S_n\cap \widetilde \Delta_n$. This fact implies that, in this case, the system has not invariant curves in $S_n$.

As we stated above, the main goal of this contribution is to prove the existence of an invariant curve when $k_1>0$ generated through a Neimark-Sacker bifurcation that occurs at the same time that both a transcritical bifurcation and the appearance of a new line of fixed points.

We remark that in \cite{Hofbauer1984} it is proved the existence of an invariant curve of amplitude $O(1/\sqrt{C})$ when $C\to \infty$. Here, instead, we are looking for an invariant curve in a different region of the space of parameters, and focus our attention on the above-mentioned functional shift.

Since we assume $C>0$, we can rewrite $F_i(x) $ as
$\frac{1+\widetilde k_ix_{i-1}}{1+\widetilde \phi(x)}x_i $
with $\widetilde \phi(x) =\sum_{i=1}^n \widetilde k_i x_ix_{i-1}$
and $\widetilde k_i= k_i/C$.
In this way we can get rid of $C$. We write $k_i$ again instead of $\widetilde k_i$.
Notice that letting $C$ go to $\infty $ in \eqref{eq63} results in letting the new parameters $k_i$ tend to zero. Here, we will only let $k_1$ go to zero keeping the other parameters fixed. Concretely, we will take $k_j>0$, for $2\le j\le n$, arbitrary and $k_1$ variable such that $k_1 > k_1^*$, with
$k_1^*= -(\sum _{j=2}^n \frac{1}{k_j})^{-1} <0 $.

\subsection{Fixed points and stability}

As a first step to understand the dynamics and for the sake of completeness, we give a brief description of the fixed points of system (\ref{eq63}) and their stability. The unique fixed point $P$ in the interior of the simplex $S_n\cap \widetilde \Delta _n$ was studied in \cite{Hofbauer1984}.
In \cite{Perona2020}  the fixed points in the boundary of the simplex were also studied.

Since the fixed points must satisfy $F_i(x)=x_i$ for all $i$, for the points in $\widetilde{\Delta}_n$, from \eqref{eq63}, we get the condition $k_ix_{i-1}=\phi(x)$,
$1\le i\le n$,  or equivalently,
\begin{equation*}
    k_2x_1=k_3x_2=\cdots=k_nx_{n-1}=k_1x_n=\phi(x).
\end{equation*}
For $k_1\le 0$, there are no fixed points in $S_n\cap \widetilde{\Delta}_n$.
If $k_1>0$, the last set of equations gives $x_i=\frac{k_1}{k_{i+1}}x_n$, $1\le i\le n-1$. Using $\sum_{i=1}^n x_i=1$, we get that the fixed point is
\begin{equation*}
P:=(p_1,\dots,p_n),\qquad \text{with } \quad p_i=\frac{1}{k_{i+1}M_1}, \quad 1\le i\le n,
\end{equation*}
where $M_1=\sum_{j=1}^n\frac{1}{k_j}$.
When $k_1=0$, $P$ coincides with $Q=(0,0\dots, 0,1)$ and when $k_1^* <k_1<0$,  $M_1<0$ and therefore $p_j<0$, $2\le j\le n$.

Moreover, in Proposition 1 of \cite{Perona2020} it was proved that
  $x\in\Delta_n\cap \widetilde{\Delta}_n$ is a fixed point if and only if $k_ix_ix_{i-1}=0$ $\forall i$. In the four species case, for $k_1\ne 0$ in the boundary of the simplex we have
the segments of fixed points $\{\,(\alpha,0,1-\alpha,0)\,| \; \alpha\in[0,1]\,\}$ and $\{\,(0, \alpha,0,1-\alpha)\, | \; \alpha\in[0,1]\,\}$. If $k_1=0$ we have the additional segment of fixed points $\{\,(\alpha,0,0,1-\alpha)\, | \; \alpha\in[0,1]\,\}$.
In particular, the vertices  $q^{(m)}:=(\delta_{m,1},\dots,\delta_{m,4})$, $1\le m \le 4$,
of the simplex $S_4$ are always fixed points.
Here $\delta_{k,l}$ is the Kronecker delta.
When $k_1\to 0^{+}$ the inner fixed point $P$ tends to the fixed point $Q=q^{(4)}=(0,0,0,1)$.

Again, in the four species case, we have that the eigenvalues of the inner fixed point $P$ are
\begin{equation*}
    \lambda_j = 1+\frac{1}{M_1+1}e^{i\theta_j}, \qquad \theta_j=e^{\frac{2\pi ij}{4}}, \qquad  1,2,3,
\end{equation*}
together with $\lambda_0= 1+\frac{1}{M_1+1}$ which has the eigenvector $(1,1,1,1)$ orthogonal to $S_4$. Therefore, concerning the dynamics in  $S_4$ the relevant eigenvalues are $\lambda_1,\,  \lambda_2$ and $\lambda_3$ (see~\cite{Hofbauer1984,Perona2020}). We also have that
\begin{equation*}
\begin{aligned}
    |\lambda_1|^2=&|\lambda_3|^2 = 1+ \Big(\frac{1}{M_1+1}\Big)^2 > 1, \\
    |\lambda_2|^2=&\Big(1- \frac{1}{M_1+1}\Big)^2<1.
\end{aligned}
\end{equation*}

Therefore, $P$ is unstable.

The eigenvalues of $q^{(m)}$ are given in \cite{Perona2020}. They are 1 (double) and $1+k_{m+1}$. In particular $Q=q^{(4)}$ has one eigenvalue that goes from bigger than 1 to less than 1 when $k_1$ goes from positive to negative.
	
\section{A non-generic Neimark-Sacker bifurcation theorem}

In this section, we recall a result by Hofbauer and Iooss  in \cite{HofbauerIooss1984}, that studies a Neimark-Sacker bifurcation for difference equations. In that paper the authors call it \emph{Hopf bifurcation} but we prefer to refer to it as \emph{Neimark-Sacker} since it seems to us that nowadays this term is more used for maps, see for instance~\cite{Kuznetsov2004}. The result deals with a discretization of a differential equation near a fixed point with two purely imaginary eigenvalues and the remaining ones with negative real part. The final goal is to prove that an invariant curve appears around the fixed point of the discrete time system.

First, we consider an autonomous differential equation
\begin{equation}\label{eq51}
    \dot{x}=f(x)
\end{equation}
defined in an open set of $\mathbb{R}^n$ and we assume that the origin is an equilibrium point, i.e. $f(0)=0$.

As in the Euler's method of integration, we consider the {following} family of maps
\begin{equation}\label{eq52}
    T_{\varepsilon}: x \mapsto x + \varepsilon f(x),\qquad \varepsilon>0 \quad  \text{small}.
\end{equation}

Since $f(0)=0$ we can write
\begin{equation}\label{eq53}
    f(x)=Ax+O(\|x\|^2),
\end{equation}
where $A=Df(0)$.
We immediately have that the maps $T_{\varepsilon}$ have the form
\begin{equation*}
    T_{\varepsilon}(x)=(\mbox{Id}+\varepsilon A)x+\varepsilon O(\|x\|^2).
\end{equation*}

It is clear that $\lambda$ is an eigenvalue of $A$ if and only if $1+\varepsilon\lambda$ is an eigenvalue of $DT_{\varepsilon}(0)$. If  $A$ has a pair of purely imaginary eigenvalues, then we have that the fixed point is unstable for the map (\ref{eq52}) for every $\varepsilon>0$, although $x=0$ could be asymptotically stable for system (\ref{eq51}).

\begin{definition}\label{def:WeaklyStable}
Assume that the origin is an equilibrium point of system \eqref{eq51} and it has a pair of purely imaginary eigenvalues {$\pm i\,\omega$}. Suppose $f$ is sufficiently differentiable and that \eqref{eq51} can be transformed, around the origin, by a change of coordinates, into the form
\begin{equation} \label{equation-nf}
    \left\{
    \begin{aligned}
    \dot{z}&=i\omega z +\sum_{j=1}^{2k}\alpha_jz|z|^{2j}+O(|z|+|v|)^{4k+2},\\
    \dot{v}&=Av+O(|z|^2+|v|^2),
    \end{aligned}\right.
\end{equation}
where $|z|^{2}=z\overline{z}$, and $\alpha_j=a_j+ib_j$.
We say that the origin is a \emph{weakly stable equilibrium point of order $k$} if there exists $k\ge 1$ such that $a_1=\dots=a_{k-1}=0$ and $a_k<0$.
\end{definition}

In \cite{HofbauerIooss1984} the following theorem is proved.

\begin{theorem}\label{th51}
Consider equation (\ref{eq51}) in an open neighbourhood of $x=0$ in $\mathbb{R}^n$ with $f$ sufficiently differentiable 
\begin{enumerate}
\item [(1)] $f(0)=0$ and $Df(0)$ has two purely imaginary eigenvalues $\pm i\omega$, and the rest of the eigenvalues have negative real part, and
\item [(2)] the equilibrium point $x=0$ is a weakly stable equilibrium point of order $k$, with $4k+2 \le r$.
\end{enumerate}
Then, for any family of maps $T_{\varepsilon}$ of class $C^r$, of the form
\begin{equation}\label{eq54}
    T_{\varepsilon}(x)=x+\varepsilon f(x)+O(\varepsilon^2\|x\|^2), \qquad  \varepsilon>0,
\end{equation}
there exists an $\varepsilon$-dependent family of invariant {and attracting} closed curves around the fixed point $x=0$ of radius
 $O({\varepsilon^{1/2k}})$.

\end{theorem}

\section{A degenerate transcritical Neimark-Sacker bifurcation}

In this section, our goal is to prove analytically the existence of an invariant curve applying Theorem \ref{th51} to our discrete-time system. In other words, we will prove that  {an invariant curve born when $k_1=0$ persists for $k_1$ positive and sufficiently small.} Since the bifurcation is very degenerate we {will} uncouple the Neimark-Sacker and the transcritical bifurcations. For {this purpose}, we force the inner fixed point to be located at the ``center" of the simplex for all values of $k_1$. This is accomplished using barycentric coordinates. Let
\begin{equation*}
    y_i=\frac{k_{i+1}x_i}{\sum_{j=1}^{4}k_{j+1}x_j}, \qquad 1\le i\le 4.
\end{equation*}
Indeed, this change allows to separate the inner fixed point {$P$} from the vertex ${Q}=(0,0,0,1)$ transforming the fixed point $p$ into the ``center point" of the simplex $S_4$:
\begin{equation*}
    p=\Big(\frac{1}{4},\frac{1}{4},\frac{1}{4},\frac{1}{4}\Big).
\end{equation*}
This transformation is singular at $k_1=0$, but it facilitates the study of the system near to the fixed point $p$, since it {brings it} far from the other fixed points.

Now, we are going to compute the new system in barycentric coordinates. First, we notice that
\begin{equation*}
    x_i=\frac{y_i\sum_{j=1}^{4}k_{j+1}x_j}{k_{i+1}}.
\end{equation*}
Since this change of coordinates sends  $\Delta_n$ to $\Delta_n$ we can write
\begin{equation*}
    \sum_{i=1}^{4}x_i=1 \quad \Longleftrightarrow \quad  \sum_{i=1}^{4}\frac{y_i\sum_{j=1}^{4}k_{j+1}x_j}{k_{i+1}}=1 \quad \Longleftrightarrow \quad \sum_{j=1}^{4}k_{j+1}x_j=\frac{1}{\sum_{i=1}^{4}\frac{y_i}{k_{i+1}}},
\end{equation*}
and we obtain
\begin{equation*}
x_i=\frac{y_i}{N(y)k_{i+1}}, \qquad \mbox{where } \quad N(y)=\sum_{j=1}^{4}\frac{y_j}{k_{j+1}}.
\end{equation*}
We can express our map $F$ in (\ref{eq63}) in the new variables $y_i$ as
\begin{equation*}
\begin{aligned}
    F_i(y)=&\frac{k_{i+1}F_i(x)}{\sum_{j=1}^4k_{j+1}F_j(x)}=\frac{k_{i+1}\Big(\frac{1+k_ix_{i-1}}{1+\phi(x)}x_i\Big)}{\sum_{j=1}^4\Big(\frac{1+k_jx_{j-1}}{1+\phi(x)}x_j\Big)k_{j+1}}=\frac{k_{i+1}\Big(1+k_i\frac{y_{i-1}}{k_iN(y)}\Big)\frac{y_i}{k_{i+1}N(y)}}{\sum_{j=1}^4\Big(1+k_j\frac{y_{j-1}}{k_jN(y)}\Big)\Big(\frac{y_j}{k_{j+1}N(y)}\Big)k_{j+1}}\\
    =& \frac{1+\frac{y_{i-1}}{N(y)}}{\sum_{j=1}^4\Big(1+\frac{y_{j-1}}{N(y)}\Big)y_j}y_i.
\end{aligned}
\end{equation*}
Next, we perform a translation to have the fixed point at the origin:
\begin{equation*}
    z_i=y_i-\frac{1}{4}, \qquad  1\leq i\leq 4.
\end{equation*}
In these new coordinates, $\sum_{j=1}^4 z_j=0$ and the value of $N(z)$ is given by
\begin{equation*}
    N(z)=\sum_{j=1}^4\frac{(z_j+\frac{1}{4})}{k_{j+1}}=\sum_{j=1}^4\frac{z_j}{k_{j+1}}+\frac{1}{4}\sum_{j=1}^4\frac{1}{k_{j+1}}   .
\end{equation*}
{Moreover, the components of $F$ become:}

\begin{equation*}
\begin{aligned}
F_i(z)&=F_i(y)-\frac{1}{4}=
\frac{1+\frac{y_{i-1}}{N(y)}}{\sum_{j=1}^4\Big(1+\frac{y_{j-1}}{N(y)}\Big)y_j}y_i-\frac{1}{4}
\frac{N(y)+y_{i-1}}{\sum_{j=1}^4\Big(N(y)+y_{j-1}\Big)y_j}y_i-\frac{1}{4}
\\
&=\frac{N(z)+z_{i-1}+\frac14}{W(z)}(z_i+\frac14)-\frac{1}{4}=z_i+\Big(\frac{N(z)+z_{i-1}+\frac14}{W(z)}-1\Big)(z_i+\frac14)
\\
&=z_i+\frac{N(z)+z_{i-1}+\frac14-W(z)}{W(z)}(z_i+\frac14),
\end{aligned}
\end{equation*}
where $W(z):=\sum_{j=1}^4\Big(N(z)+z_{j-1}+\frac14\Big)(z_j+\frac14)$. Note that $$
W(z)=\sum_{j=1}^4 N(z) (z_j+\frac14)+\sum_{j=1}^4 (z_{j-1}+\frac14)(z_j+\frac14)=N(z)+\sum_{j=1}^4 z_{j-1}z_j+\frac14.
$$
Therefore,
\begin{equation*}
F_i(z)=z_i+\frac{z_{i-1}-\sum_{j=1}^4z_jz_{j-1}}{\frac{1}{4}\Big(1+\sum_{j=1}^4\frac{1}{k_{j+1}}\Big)+\sum_{j=1}^4\frac{z_j}{k_{j+1}}+\sum_{j=1}^4z_{j-1}z_j}\left(z_i+\frac{1}{4}\right).
\end{equation*}
Now, keeping the lower order terms, we can rewrite the {components of the} system as:
\begin{equation}\label{eqFiz}
\begin{array}{rl}
    F_i(z)&=z_i+\delta(z_{i-1}-\sum_{j=1}^4 z_jz_{j-1})\dfrac{1/4+z_i}{1/4+z_4}+O(\delta^2)O(|z|^2)\\
          &=z_i+\delta(z_{i-1}-\sum_{j=1}^4 z_jz_{j-1})(1+4\,z_i)(1-4\,z_4+16\,z_4^2+O(z_4^3))+O(\delta^2)O(|z|^2),\\
\end{array}
\end{equation}
where
$$\delta=\frac{k_1}{1+k_1(1+M_2)}, \qquad \text{with } \quad {M_2}=\frac{1}{k_2}+\frac{1}{k_3}+\frac{1}{k_4}.$$
Note that $O(\delta)=O(k_1)$ and, more importantly, that $F$ expands exactly as $T_{\varepsilon}$ in \eqref{eq54} of Theorem \ref{th51}. In order to apply this theorem, we first reduce the dimension of the map $F$ by $1$ using that $\sum_{j=1}^4 z_j=0$. We choose to eliminate $z_2$; as a consequence, we also have that
$$\sum_{j=1}^4z_jz_{j-1}=(z_1+z_3)(z_2+z_4)=-(z_1+z_3)^2.$$
Let us call $G$ the new map, so that $G_i(z_1,z_3,z_4)=F_j(z_1,-z_1-z_3-z_4,z_3,z_4)$, with
$G_1=F_1$, $G_2=F_3$ and $G_3=F_4$.
Then, we can express the system as a family of maps $G(z)$ of the form
\begin{equation*}
    G(z)=z+\delta g(z)+O(\delta^2)O(|z|^2),
\end{equation*}
where the components of $g$ are obtained from \eqref{eqFiz} by substituting $z_2=-z_1-z_3-z_4$:
\begin{equation}\label{eq64}
    \left\{
    \begin{aligned}
    g_1(z)&=(1+4z_1)(z_4+(z_1+z_3)^2)(1-4z_4+16z_4^2+O(z_4^3)),\\
    g_2(z)&=(1+4z_3)(-z_1-z_3-z_4+(z_1+z_3)^2)(1-4z_4+16z_4^2+O(z_4^3)),\\
    g_3(z)&=z_3+(z_1+z_3)^2+O(z_4^3).\\
    \end{aligned}
    \right.
\end{equation}

Expanding in powers of $z$ we can write 
\begin{equation}\label{eq64b}
    \left\{
    \begin{aligned}
    g_1(z)&=z_4+P_{12}(z)+P_{13}(z)+O(|z|^4),\\
    g_2(z)&=(-z_1-z_3-z_4)+P_{22}(z)+P_{23}(z)+O(|z|^4),\\
    g_3(z)&=z_3+P_{32}(z)+P_{33}(z)+O(|z|^4),\\
    \end{aligned}
    \right.
\end{equation}
where $P_{ij}$ indicates the term of degree $j$ in the $i$-th component of the vector field, and
$$
\begin{array}{l}
P_{12}(z)=-4z_4^2+4z_4z_1+(z_1+z_3)^2, \\ 
P_{13}(z)=4\,(z_1-z_4)\,(z_1 + z_3 + 2\,z_4)\,(z_1 + z_3 - 2\,z_4),\\
P_{22}(z)=4z_4(z_1+z_3+z_4)-4 z_3 (z_1+z_3+z_4)+(z_1+z_3)^2, \\
P_{23}(z)=4\,(z_3-z_4)\,(z_1 + z_3 + 2\,z_4)^2,\\
P_{32}(z)=(z_1+z_3)^2,   \qquad P_{33}(z)= 0.
\end{array}
$$

Now, we have to {check that $g$ satisfies} the two hypotheses of Theorem \ref{th51}. Clearly, we have that {both} $g(0)=0$ and the derivative of $g$ at the origin,
\begin{equation*}
    Dg(0)=\left(\begin{array}{ccc}
     0 & 0 & 1 \\
    -1 & -1 & -1 \\
     0 & 1 & 0
    \end{array}\right),
\end{equation*}
{has two purely imaginary eigenvalues and a third one whose real part is negative; more precisely, the eigenvalues are $\lambda_{1,2}=\pm i$, and $\lambda_3 = -1$. Thus, the first hypothesis of the theorem follows.}
The corresponding eigenvectors are $v_1=(1,-1,i)$, $v_2=(1,-1,-i)$ and $v_3=(1, \hspace{3mm}1,-1).$

For the second hypothesis of the theorem, we need to compute the normal form for the system $\dot{z}=g(z)$. First, we diagonalize $Dg(0)$ using the linear change $z=C\,\zeta$, where 
\begin{equation*}
    C=\left(\begin{array}{ccc}
     1 & 1 & 1 \\
    -1 & -1 & 1 \\
     i & -i & -1
    \end{array}\right) \mbox{\quad and \quad}
    C^{-1}=
    \frac{1}{4}
    \left(\begin{array}{ccc}
     1-i & -1-i & -2i \\
    1+i & -1+i & 2i \\
     2 & 2 & 0
    \end{array}\right).
\end{equation*}
{In the new set of variables $\zeta=(\xi,\overline{\xi},\eta)$, system $\dot{z}=g(z)$ is transformed into}
\begin{equation*}
    \dot{\zeta}=g^{(1)}(\zeta):=C^{-1}g(C\zeta).
\end{equation*}
Observe that, by construction, the linear term of $g^{(1)}$ becomes
$\left(\begin{array}{c}
      i\xi \\
      -i\overline{\xi} \\
      -\eta
    \end{array}\right).$
    
Next, we have to compute $g^{(1)}$ for quadratic and cubic terms. We first compute the corresponding terms in $g(C\zeta)$; writing each component $i$ in the form
$g_{i2}(C\zeta)+g_{i3}(C\zeta)$, we have
$$
\begin{aligned}
    g_{12}(C\zeta)&=4(-\eta^2+ \xi^2(1+i)-2\xi\overline{\xi}+\xi\eta(-1+3i)+\overline{\xi}\eta(-1-3i)+\overline{\xi}^2(1-i)),\\
    g_{13}(C\zeta)&=
16\,i\,(\xi\,i - \overline{\xi}\,i - 2\,\eta - \xi - \overline{\xi})\,(\xi\,i - \overline{\xi}\,i - 2\,\eta)\,(\xi - \overline{\xi}),\\
    g_{22}(C\zeta)&=4(-\eta^2+ \xi^2(-1+i)+2\xi\overline{\xi}+\xi\eta(1-i)+\overline{\xi}\eta(1+i)-\overline{\xi}^2(1+i)),\\
    g_{23}(C\zeta)&=16\,(\xi\,i - \overline{\xi}\,i - 2\,\eta + \xi + \overline{\xi})\,(\xi - \overline{\xi})^2,\\
g_{32}(C\zeta)&=4\eta^2,\\
g_{33}(C\zeta)& = 0.
\end{aligned}
$$

Once we have $g(C\zeta)$, we then compute $g^{(1)}(\zeta)=C^{-1}g(C\zeta)$:
\begin{equation*}
    g^{(1)}(\zeta)=\frac{1}{4}\left(\begin{array}{c}
      (1-i)g_1(C\zeta)+(-1-i)g_2(C\zeta)-2ig_3(C\zeta) \\
      (1+i)g_1(C\zeta)+(-1+i)g_2(C\zeta)+2ig_3(C\zeta)\\
      2g_1(C\zeta)+2g_2(C\zeta)
    \end{array}\right)=:
\left(\begin{array}{c}
g^{(1)}_1(\zeta)\\
g^{(1)}_2(\zeta)\\
g^{(1)}_3(\zeta)\\
    \end{array}\right).
\end{equation*}
We decompose the three components as $g^{(1)}_i=g^{(1)}_{i1}+g^{(1)}_{i2}+g^{(1)}_{i3}$, for $i=1,2,3$. Clearly $g^{(1)}_{11}(\zeta)=i\xi$, $g^{(1)}_{21}(\zeta)=-i\overline{\xi}$ and $g^{(1)}_{31}(\zeta)=-\eta$, and
$$
\begin{aligned}
    g^{(1)}_{12}(\zeta)=&4\left(\xi^2-\xi\overline{\xi}+i\xi\eta-(1+i)\overline{\xi}\eta\right),\\
    g^{(1)}_{13}(\zeta)=&16\left(-i\xi^3+2i\xi^2\overline{\xi}+2\xi^2\eta-i\xi\overline{\xi}^2+(-3+i)\xi\overline{\xi}\eta+(1+i)\xi\eta^2+(1-i)\overline{\xi}^2\eta\right.\\&\qquad\left.-(1+i)\overline{\xi}\eta^2\right),\\
    g^{(1)}_{22}(\zeta)=&4\left(-\xi\overline{\xi}+(-1+i)\xi\eta+\overline{\xi}^2-i\overline{\xi}\eta\right),\\
    g^{(1)}_{23}(\zeta)=&16\left(i\overline{\xi}^3+i\xi^2\overline{\xi}+(1+i)\xi^2\eta-2i\xi\overline{\xi}^2-(3+i)\xi\overline{\xi}\eta+(-1+i)\xi\eta^2+2\overline{\xi}^2\eta\right.\\&\qquad\left.+(1-i)\overline{\xi}\eta^2\right),\\
    g^{(1)}_{32}(\zeta)=&4\left(i\xi^2+i\xi\eta-i\overline{\xi}^2-i\overline{\xi}\eta-\eta^2\right),\\
    g^{(1)}_{33}(\zeta)=&16\left(\xi^3-\xi^2\overline{\xi}+(1+i)\xi^2\eta-\xi\overline{\xi}^2-2\xi\overline{\xi}\eta+2i\xi\eta^2+\overline{\xi}^3+(1-i)\overline{\xi}^2\eta-2i\overline{\xi}\eta^2\right).
\end{aligned}
$$

We now proceed to compute the normal form of $g^{(1)}(\zeta)$ by means of a generic change of coordinates of quadratic order that kills all quadratic terms (which are non-resonant) of $g^{(1)}(\zeta)$ and preserves the linear ones. Let
\begin{equation*}
    h({\rm x})={\rm x}+\widetilde{h}({\rm x}),
\end{equation*}
where
\begin{equation*}
    \widetilde{h}({\rm x})=\left(\begin{array}{c}
      a_{200}x^2+a_{020}y^2+a_{002}z^2+a_{110}xy+a_{101}xz+a_{011}yz \\
      b_{200}x^2+b_{020}y^2+b_{002}z^2+b_{110}xy+b_{101}xz+b_{011}yz\\
      c_{200}x^2+c_{020}y^2+c_{002}z^2+c_{110}xy+c_{101}xz+c_{011}yz
    \end{array}\right),
\end{equation*}
and consider the change of variables $\zeta=h({\rm x})$, with ${\rm x}=(x,y,z)$. We have that
\begin{equation*}
    \dot{\zeta}=Dh({\rm x})\dot{{\rm x}},
\end{equation*}
and so
\begin{equation}\label{UltimCanvi}
    \dot{{\rm x}}=Dh({\rm x})^{-1}g^{(1)}(h({\rm x}))=:g^{(2)}({\rm x}).
\end{equation}

\begin{remark}\label{rem01}
To do the computations, we only keep track the terms up to degree $3$ and we take advantage of the degree-structure presented in the previous steps to discard terms of degree $4$ or higher in $g^{(1)}(h({\rm x}))$. Moreover, we can approximate $Dh({\rm x})^{-1}$ by
\begin{equation*}
    Dh({\rm x})^{-1}\approx \Id-D\widetilde{h}({\rm x})+D\widetilde{h}({\rm x})^2.
\end{equation*}
Note that, when we substitute this approximation in \eqref{UltimCanvi}, the $\Id$ applies to the expression of $g^{(1)}(h({\rm x}))$ up to degree $3$, but $-D\widetilde{h}({\rm x})$ applies only up to quadratic terms and $D\widetilde{h}({\rm x})^2$ only to the linear terms. 
\end{remark}

Following the strategy commented in Remark \ref{rem01}, the quadratic terms of the new system \eqref{UltimCanvi} have the following components
\begin{equation*}
\begin{aligned}
    g^{(2)}_{1}({\rm x})=&(4-ia_{200})x^2+3ia_{020}\,y^2+(ia_{002}+2a_{002})z^2+(ia_{110}-4)xy\\
    &+(4i+a_{101})xz+(-4-4i+2ia_{011}+a_{011})yz,\\
    g^{(2)}_{2}({\rm x})=&(-3ib_{200})x^2+(ib_{020}+4)y^2+(-ib_{002}+2b_{002})z^2+(-ib_{110}-4)xy\\
    &+(-2ib_{101}-4+4i+b_{101})xz+(-4i+b_{011})yz,\\
    g^{(2)}_{3}({\rm x})=&(-2ic_{200}+4i-c_{200})x^2+(2ic_{020}-4i-c_{020})y^2+(c_{002}-4)z^2+(-c_{110})xy\\
    &+(-ic_{101}+4i)xz+(ic_{011}-4i)yz.
\end{aligned}
\end{equation*}
In order to kill every quadratic term, we must take
{\small{
\begin{equation*}
    \begin{aligned}
        &a_{200}= -4i, &\hspace{0mm} &a_{020}= 0, &\hspace{0mm}   & a_{002}= 0, &\hspace{0mm} &a_{110}= -4i, &\hspace{0mm}
        &a_{101}= -4i, &\hspace{0mm} &a_{011}= \frac{4}{5}(3-i),\\
        &b_{200}= 0, &\hspace{0mm} &b_{020}= 4i, &\hspace{0mm}
        &b_{002}= 0, &\hspace{0mm} &b_{110}= 4i, &\hspace{0mm}
        &b_{101}= \frac{4}{5}(3+i), &\hspace{0mm} &b_{011}= 4i,\\
        &c_{200}= \frac{4}{5}(2+i), &\hspace{0mm} &c_{020}= \frac{4}{5}(2-i), &\hspace{0mm}
        &c_{002}= 4, &\hspace{0mm} &c_{110}= 0, &\hspace{0mm}
        &c_{101}= 4, &\hspace{0mm} &c_{011}= 4.
\end{aligned}
\end{equation*}
}}

\noindent Next, we substitute the above values of the coefficients of $\tilde{h}({\rm x})$ into the cubic terms of ${g}^{(2)}$. It is worth mentioning that, in principle, this cubic terms can have non-zero coefficients for all the monomials. Thus, in order to have the cubic normal form, we should continue with a new change of variables that would kill all cubic terms but the resonant ones. However, by the normal form theory, we know that this new change would keep all resonant terms invariant. Since we are only interested in the sign of the real part of one specific resonant term, we do not need to perform the full change of variables. Therefore, if we call $(z,\bar{z},\nu)$ the new set of variables, we can assert that the system writes as
\begin{equation}\label{equation-nf0}
    \left\{
    \begin{aligned}
        \dot{z}&=iz+\Big(-\frac{16}{5}-\frac{48}{5}i\Big)z^2\bar{z}+\dots ,\\
        \dot{\bar{z}}&=-i\bar{z}+\Big(-\frac{16}{5}+\frac{48}{5}i\Big)z\bar{z}^2+\dots  ,\\
        \dot{\nu}&=-\nu+\frac{64}{5}z\bar{z}w+\dots .
    \end{aligned}
    \right.
\end{equation}

Observe that \eqref{equation-nf0} corresponds to the normal form \eqref{equation-nf} with $n=3$, $A=-1$ and, most importantly, $\alpha_1=-\frac{16}{5}-\frac{48}{5}i$. Since $\mbox{Re}(\alpha_1)$ is negative, from Definition \ref{def:WeaklyStable} we can ensure that the origin is a weakly stable equilibrium point of order $1$ and so we have checked the second hypothesis of Theorem \ref{th51} for our system. Therefore, we conclude that the four-dimensional discrete-time hypercycle (\ref{eq63}) presents a family of attracting invariant curves depending on the parameter $\varepsilon=k_1$, when $k_1>0$. Going back to the original variables we have that, for $k_1>0$ small, the system has a closed invariant curve which arises from $Q=q^{(4)}=(0,0,0,1)$, while, at the bifurcation value $k_1=0$, a line of fixed points appears and this corner point $Q=q^{(4)}$ collides with the inner fixed point $P$ in a transcritical bifurcation.

\section{Conclusions}

The main goal of this work was to provide an analytical proof of the existence of an attracting invariant curve in the four member discrete-time hypercycle when a cooperation coefficient approaches the functional shift, motivated by numerical evidences that were described in \cite{Perona2020}. In the discrete-time hypercycle model, this phenomenon is reflected in the fact that the parameter $k_1$ goes from positive to negative: the invariant curve shrinks to a corner of the domain and disappears throughout a Neimark-Sacker bifurcation when $k_1=0$. 

We have studied analytically this degenerate bifurcation. For this purpose, we have followed a result by Hofbauer and Iooss~\cite{HofbauerIooss1984} that provides sufficient conditions for a Neimark-Sacker bifurcation. In fact, the theorem by Hofbauer and Iooss was introduced to prove the existence of another invariant curve in the same model. However, the application of this theorem is not straightforward for the case $k_1=0$.

The coincidence of the Neimark-Sacker bifurcation with a transcritical one forced us to decouple them. For this purpose, we performed a singular change of coordinates that ensured a constant distance between the fixed points that are relevant in each bifurcation. Subsequently, in order to prove the hypotheses of the theorem by Hofbauer and Iooss, we brought the system to its normal form by making a new change of variables that eliminates all quadratic terms and reveals the resonant cubic term. By undertaking this analytical exploration, we have been able to provide a complete understanding of how the invariant curve arises in the scenario of transition from cooperation to degradation. The presence of invariant attracting curves ensures the survival of all species; the dynamics within these invariant curves is an interesting continuation of this problem that would shed light on how oscillations in the model are structured.  


\section*{Acknowledgments}
EF has been funded by the Spanish grant PID2021-125535NB-I00 (MICINN/FEDER,UE). AG has been unded by MCIN/AEI/10.13039/501100011033 and by ERDF "A way of making Europe" grants PID-2021-122954NB-I00 and PID2022-137708NB-I00, and the AGAUR project 2021SGR1039. JS has been supported by the Ram\'{o}n y Cajal grant RYC-2017-22243 funded by MCIN/AEI/10.13039/501100011033 ”FSE invests in your future”, and by the 2020-2021 Biodiversa+ and Water JPI joint call under the BiodivRestore ERA-NET Cofund (GA N°101003777) project MPA4Sustainability with funding organizations: Innovation Fund Denmark (IFD), Agence Nationale de la Recherche (ANR), Fundaçao para a Ciencia e a Tecnologia (FCT), Swedish Environmental Protection Agency (SEPA), and grant PCI2022-132926 funded by MCIN/AEI/10.13039/501100011033 and by the European Union NextGeneration EU/PRTR. This work has been also funded through the Severo Ochoa and Mar\'ia de Maeztu Program for Centers and Units of Excellence in R\&D (CEX2020-001084-M). We thank CERCA Programme/Generalitat de Catalunya for institutional support.

\end{document}